\documentclass[12pt]{article}
\font\sdopp=msbm10
\def\CI {\sdopp {\hbox{C}}}
\title{Renormalizing iterated 
elementary \\ mappings and 
correspondences of $\hbox{\tt\Large C}^2$
}
\author{Claudio Meneghini
\footnote{\tt The author's e-mail address:
clamen@dimat.unipv.it}
}
\begin{document}
\maketitle
\bibliographystyle{plain} 
\parindent=8pt
%caratteri cirillici
\font\cir=wncyb10
\def\Iu{\cir\hbox{YU}}
\def\Ze{\cir\hbox{Z}}
\def\pe{\cir\hbox{P}}
\def\Ef{\cir\hbox{F}}
\def\CIRC{\mathop{\tt o}\limits}
\def\quan{\vrule height6pt width6pt depth0pt}
\def\QUAN{\ \quan}
\def\BETA{\mathop{\beta}\limits}
\def\GAMMA{\mathop{\gamma}\limits}
\def\VI{\mathop{v}\limits}
\def\UI{\mathop{u}\limits}
\def\VII{\mathop{V}\limits}
\def\WI{\mathop{w}\limits}
\def\ZETA{\mathop{Z}\limits}

\newtheorem{definition}{Definition}%[chapter]
\newtheorem{defi}[definition]{D\'efinition}%[chapter]
\newtheorem{lemma}[definition]{Lemma}
\newtheorem{lemme}[definition]{Lemme}
\newtheorem{proposition}[definition]{Proposition}
\newtheorem{theorem}[definition]{Theorem}        
\newtheorem{theoreme}[definition]{Th\'eor\`eme}        
\newtheorem{corollary}[definition]{Corollary}  
\newtheorem{corollaire}[definition]{Corollaire}  
\newtheorem{remark}[definition]{Remark}  
\newtheorem{remarque}[definition]{Remarque}
  
\font\sdopp=msbm10
\def\ERRE {\sdopp {\hbox{R}}}
\def\QU {\sdopp {\hbox{Q}}}
\def\CI {\sdopp {\hbox{C}}}
\def\DI {\sdopp {\hbox{D}}}
\def\ENNE{\sdopp {\hbox{N}}}
\def\ZETA{\sdopp {\hbox{Z}}}
\def\PI {\sdopp {\hbox{P}}}
\def\M{\hbox{\tt\large M}}
\def\N{\hbox{\boldmath{}$N$\unboldmath}} 
\def\P{\hbox{\boldmath{}$P$\unboldmath}} 
\def\tr{\hbox{\boldmath{}$tr$\unboldmath}} 
\def\f{\hbox{\large\tt f}} 
\def\F{\hbox{\boldmath{}$F$\unboldmath}} 
\def\G{\hbox{\boldmath{}$G$\unboldmath}} 
\def\L{\hbox{\boldmath{}$L$\unboldmath}} 
\def\h{\hbox{\boldmath{}$h$\unboldmath}} 
\def\e{\hbox{\boldmath{}$e$\unboldmath}} 
\def\g{\hbox{\boldmath{}$g$\unboldmath}} 
\def\u{\hbox{\boldmath{}$u$\unboldmath}} 
\def\v{\hbox{\boldmath{}$v$\unboldmath}} 
\def\U{\hbox{\boldmath{}$U$\unboldmath}} 
\def\V{\hbox{\boldmath{}$V$\unboldmath}} 
\def\W{\hbox{\boldmath{}$W$\unboldmath}} 
\def\id{\hbox{\boldmath{}$id$\unboldmath}} 
\def\alph{\hbox{\boldmath{}$\alpha$\unboldmath}} 
\def\bet{\hbox{\boldmath{}$\beta$\unboldmath}} 
\def\gam{\hbox{\boldmath{}$\gamma$\unboldmath}} 
\def\pphi{\hbox{\boldmath{}$\varphi$\unboldmath}} 
\def\ppsi{\hbox{\boldmath{}$\psi$\unboldmath}} 
\def\Ppsi{\hbox{\boldmath{}$\Psi$\unboldmath}} 
\def\brevve{}
\def\labelle #1{\label{#1}}

\begin{abstract}
After proving a multi-dimensional extension
of Zalcman's renormalization lemma and considering maximality problems about dimensions,
we find renormalizing polynomial families 
for iterated elementary mappings,
extending this result to some kinds of correspondences (by means of 'algebraic' 
renormalizing families) and to the family of the iterated mappings of an automorphism of $\CI^2$ admitting a repulsive fixed point (by means of a family of polunomial automorphisms composed with a Fatou-Bieberbach one).
All families will allow maximal- dimension renormalizations.
\end{abstract}
\section{\bf \large Foreword} 
Renormalizing holomorphic families
is a useful tool in complex analysis: 
in one variable, for instance, it allows, via
'Zalcman's renormalization lemma' (see \cite{zalcman}, \cite {berteloot}, pag. 9
or lemma
\ref{zalcman} for 
a more general version) a quite direct proof of both great Picard's and Montel's theorems. 
The process of renormalization could be described as 
composing on the right the elements
of a nonnormal family 
${\cal F}$
with a family of  contracting mappings 
(depending on the nature of the problem: affine, polynomial, biholomorphic ones etc...) and then
extracting a normally convergent subsequence from 
the new family.
It is worthwhile to recall that renormalization via
Zalcman's lemma
is made up by means of affine functions and
yields an {\it entire} limit function: this is
crucial in proving the quoted theorems (see again \cite{berteloot}).
Moreover, as shown in \cite{milnor}, chap. 8,
it is always possible
to renormalize {\it linearly} the
familiy of iterates
of an endomorphism of $\PI^1$
admitting a repulsive fixed point:
however, in general, this will not 
yield an {\it entire}
function.

The situation is different in more than one variable: the following example, adapted from \cite{rudin} (9.2) shows, for instance, the existence of nonnormal families
of iterates
of
automorphisms of $\CI^2$, admitting a repulsive fixed point in $0$,
which are by no means linearly renormalizable.

Let $F\in Aut(\CI^2)$ be defined by
$
F(z,w)=(\alpha z, \beta w+z^2)
$,
with $\vert\alpha\vert>1$, $\vert\beta\vert>1$;
$F$ admits a repulsive fixed point in $0$, hence $\{F^{\CIRC k}\}$ cannot be normal.

Now
$
F^{\CIRC k}(z,w)=
(\alpha^k z, \beta^k w+
\beta^{k-1}[1+c+\cdots c^{k-1}]z^2)
$,
where $c=\alpha^2/\beta$.
Suppose 
$\vert\beta\vert<
\vert\alpha\vert^2$: then
$[F_*(0)]^{-k}F^{\CIRC k}(z,w)=
(z,w+
\beta^{-1}[1+c+\cdots c^{k-1}]z^2)$;
since $c>1$, the coefficient of $z^2$
in $([F_*(0)]^{-k}F^{\CIRC k})_2$
diverges.
Thus, 'dividing' by the differentials in $0$ allows in general no kind of renormalization.
\vskip0,1truecm
By contrast, Zalcman's renormalization lemma holds in higher dimension too (lemma
\ref{zalcman}); but, at least as far as
{\it affine renormalization} is concerned,
once a renormalized family has been got, its limit function will generally have
lower-dimensional image.
This fact seems to suggest the existence of (complex) directions with a 'major degree'
of nonnormal behaviour: this is not quite satisfactory, since one should also compare different directions.

By restricting our investigations to holomorphic dynamics, we shall show that
a nonnormal family ${\cal I}$ of 
iterates 
of a so-called
{\it elementary mapping} $E:\CI^2\rightarrow\CI^2$ can be globally renormalized by means of a family
${\cal P}$ of polynomial mappings whose degrees are uniformly bounded by a quantity depending only on the eignevalues of $E_*$, achieving in fact a limit mapping whose image is two dimensional: this allows us to 
gain information about the chaotic behaviour if ${\cal I}$ by comparison
with the contracting one of ${\cal P}$,  easier to study.

Let $G$ be the group of all polynomial automorphisms of $\CI^2$:
we recall that an element of $G$
is conjugated, within this group, to an elementary automorphism or to a H\'enon mapping, according to its dynamical degree being bounded or not (see \cite{friedland}):
in spite of being
the dynamics of {\sf polynomial} elementary mappings well known, our 
techniques
still keep some interest since they are
valid for nonpolynomial and for some kind of 'multi-valued' ones too (see section \ref{main}). Moreover, we shall show, as an application of our main result, that the family of the iterates
of an automorphism $H$ 
of $\CI^2$  
(hence, in particular, a H\'enon mapping)
with a repulsive fixed point $p$
admits a globally renormalizing family of mappings on the {\sf repelling basin} of $p$, i.e. the attracting basin of $p$ with respect to $H^{-1}$.
\section{\bf\large Zalcman's renormalization lemma extended }
\begin{lemma} 
$
\hbox{\tt Metric-space lemma
(see \rm\cite{gromov}, \tt pag 256
and \rm \cite{berteloot2}).\it\ Let
}$
\\
$(X,d)$ be a complete metric space and $M:X\rightarrow \ERRE^+_0$ 
a locally bounded 
%\hbox{
function: then 
for each $\sigma>0$ and 
$u\in M^{-1}(\ERRE^+)$
there
exists
$v\in X$
such that:
%} 
\\
{\tt (i)} $d(u,v)\leq 2/(\sigma M(u))$,
{\tt (ii)} $M(v)\geq M(u)$;\\
{\tt (iii)} $d(x,v)\leq 1/(\sigma M(v))\Rightarrow
M(x)\leq 2M(v)$.
\labelle{metric}
\end{lemma}
{\bf Proof:} were there no such $v$, we could make up a sequence $\{v_n\}$ such that $v_0=u$, $M(v_{n+1})\geq 2M(v_n)\geq
2^{n+1}M(u)$ and $d(v_{n+1},v_n)\leq
(\sigma M(v_n))^{-1}\leq 2^{-n}(\sigma M(u))^{-1}$, implying that this sequence is Cauchy like: this is a contradiction.
\QUAN
\begin{lemma}
A family of 
${\cal F}:=\{f_{\alpha}\}\subset
{\cal O}(\DI^N, \PI^N)$
{\sf is not} normal at $v\in\DI^N$ if and only if
there exist  
sequences $\{v_n\}\subset\DI^N$ $\{v_n\}\rightarrow v$,
$\{r_n\}\subset\ERRE^+$, with $\{r_n\}\rightarrow 0$,
and $\{f_n\}\subset{\cal F}$ such that
$\{f_n(v_n+r_nw)\}$ converges normally to a nonconstant entire mapping $h$ such that,
for each $w\in\CI^N$,
$\vert h_*(w)\vert
\leq 2$, with $\vert h_*(w)\vert=
\max_{i=1\cdots N}[\g(\partial h/\partial u_i, \partial h/\partial u_i)]^{1/2}$, where $\g$ is 
Fubini-Study's metric on $\PI^N$.
\labelle{zalcman}
\end{lemma}
{\bf Proof:} 
we prove only $(\Rightarrow)$.
We can find sequences
 $\{\xi_n\}\rightarrow v$ 
in $\DI^N$ and 
$\{f_n\}\subset{\cal F}$ 
such that $\vert (f_n)_*(\xi_n)
\vert\geq n^2$.
We may suppose that $\{\xi_n\}$
is contained in a closed subset $X$ of $\DI^N$.
For each $n$, apply lemma \ref{metric} to $X$
with the euclidean metric, $M(x)=\vert 
(f_n)_*(x)\vert$, $u=\xi_n$
and $\sigma=1/n$: an element $v_n\in X$ is yielded such that:
{\tt (a)} $d(\xi_n,v_n)\leq 2/n$
and 
{\tt (b)} $\vert 
(f_n)_*(x)\vert\leq 2$ if 
$d(x,v_n)\leq n/\vert 
(f_n)_*(v_n)\vert$.
Now set $r_n=\vert 
(f_n)_*(v_n)\vert
^{-1}$: then, if $h_n(w)=f_n(v_n +r_n w)$,
the family $\{g_n\}$ 
is normal,  since by {\tt (b)} 
$\vert(h_n)_*\vert\leq 2$ on $B(0,n)$.
Thus we can extract from 
$\{h_n\}$ a normally convergent
 family, whose limit we shall
 call $h$.  By {\tt (a)}
$v_n\rightarrow v$ and, 
by construction,
$\vert{h_n*}(0)
\vert\equiv 1$, 
hence $\vert {h_{n*}}(0)\vert=1$,
proving that $h$ is not a 
constant function.
Finally, $\vert {h_{n*}}(U)
\vert =\lim\vert h_{n*}(U)
\vert\leq 2$:
this ends the proof.
\QUAN
\vskip0,1truecm
{\bf Remark:} consider the example in the foreword section: by nonnormality at $0$ and lemma \ref{zalcman} there exist: $\{(z_n,w_n)\}\rightarrow 0$, $\{\gamma_n\}\rightarrow 0$ in $\ERRE^+$
such that 
\begin{equation}
(u,v)
\mapsto
\left\{
\left(\alpha^n
(z_n+\gamma_n u),
\beta_n(w_n+\gamma_n v)+
\beta^{n-1}(\sum_{l=0}^{n-1}c^l)
(z_n+\gamma_n u)^2
\right)
\right\}
\labelle{aster}
\end{equation}
converges
normally
to an entire function $G$ on $\CI^2$.
However, since the hypothesis
$\vert\beta\vert<\vert\alpha\vert^2$ implies
$\beta^{n-1}\sum_{l=0}^{n-1}c^l\sim \alpha^{2n}$ as $n\to\infty$, it is easily seen that a necessary condition for convergence in
(\ref{aster}) is
that
$\vert\gamma_n\vert$
decreases as fast as $\min(
\vert\alpha\vert^{-n},
\vert\beta\vert^{-n})$.
This easily implies 
dim($G(\CI^2)$)=1, unless
$\vert\alpha\vert=
\vert\beta\vert$.
This difficulty will be overcome
by theorem \ref{principal}.
\section{\bf\large Some definitions}
We recall that $G=(g_1,g_2)$ is an {\sf elementary mapping} of  $\CI^2$ if 
$g_1(z)=c_1 z_1$, 
$g_2(z)=c_2 z_2+h(z_1)$,
where the $c_k$'s are complex constants and
$h$ is a holomorphic 
function of 
$z_1$; $G$ is an automorphism if and only if each $c_k$
is nonzero.
A {\sf germ of elementary mapping}
will be a germ of mapping in $\CI^2$ of the form
$\G=\L+(0,\h)$, where $\L$ is the germ at $0$
of the linear mapping
$L_1(z)=c_1 z_1$, 
$L_2(z)=c_2 z_2$ and
$\h$ is a germ of holomorphic function, depending only on $u$ and with $\h(0)=0$.
{\sf The branched maximal analytical continuation of $\G$} is the holomorphic mapping $G$ defined on $S\times\CI$ (where
$\left(S,\pi   \right)$ is the branched maximal analytical continuation of $\h$; see e.g.
\cite{cassa}, chap.5) by setting $G(x,y)=L(\pi(x),y)+(0,h(x))$.
The germ $\G$ will
be told {\sf complete} if so is $\h$, i.e. if $\pi(S)=\CI$.
\begin{definition}
{\sf An elementary correspondence} 
${\cal G}$
of $\CI^2$
is the correspondence generated by a complete germ 
$\G$
of an elementary mapping, that is to say 
$(x_1,x_2){\cal R}(y_1,y_2)$ if and only if 
there exists $\widetilde x_1\in S$ such that
$x_1=\pi(\widetilde x_1)$, and 
$G(\widetilde x_1,x_2)=L(x_1,x_2)+
(0,h(\widetilde x_1))$, where $G$ is the branched maximal analytical continuation of $\G$.
\labelle{elemmap}
\end{definition}

If ${\cal R}$ is a correspondence between the sets 
$X$ and $Y$, ${\cal S}$ between 
$Y$ and $Z$ then we can define their 
{\sf composition}: if $x\in X$ and $z\in Z$, then
$x\,({\cal R}\CIRC {\cal S}) \, y$ if and only if there exists $y\in Y$ such that  
$x\, {\cal R}\, y$ and $y\, {\cal S}\, z$.
Note that it well could be ${\cal R}\CIRC {\cal S}=\emptyset$: however,
an elementary correspondence can be always iterated
preserving completeness and
without yielding the empty set.
\begin{definition}
\label{fixbra}
A correspondence ${\cal R}$ on a set $X$ {\sf has a fixed point} at
$x\in X$ if $x\,{\cal R}\, x$; 
it {\sf possesses a branch} $\phi$
in a superset ${\cal U}$ of $\{x\}\subset X$ if
$\phi $
is a {\sf function} on ${\cal U}$
such that for every $\xi\in{\cal U}$, $(\phi(\xi)=\eta)\Rightarrow \phi{\cal R}\eta$ 
\end{definition}
\section{\bf\large
Elementary mappings
\labelle{main}\labelle{mapm}
}
\begin{lemma}
If $F:\CI^2\rightarrow\CI^2$ is defined by setting $F(u,v)= (\alpha u,\beta v+h(u))$, where $h(u)=\sum_{l=0}^{\infty}\eta_l u^l$ and $\alpha$, $\beta$ are complex constants, then,
for $n\geq 0$, 
$F^{\circ n}(u,v)=
\left(
\alpha^n u,
\left[
\sum_{k=0}^{n-1}
\beta^k
h(\alpha^{n-1-k} u)
\right]+\beta^n v
   \right)$.
\vskip0,1truecm
{\bf Proof:} \rm by induction on $n$.
\QUAN
\labelle{elemcomp}
\end{lemma}

As a consequence, if $\alpha\not=0$, $\beta\not=0$, then $F$ is invertible and
there holds
$F^{-n}(u,v)=
(
\alpha^{-n} u,
[
-
\sum_{k=1}^{n}
\beta^{-k}
h(\alpha^{-n-1+k} u)
]+\beta^{-n}v
)
$.

Let now $P_N(u)=
\sum_{l=0}^{N}\eta_l u^l$
be the $N$-th degree truncation of the development of $h$
and 
$F_N^{-n}(u,v)=
(
\alpha^{-n} u,
[
-
\sum_{k=1}^{n}
\beta^{-k}
P_N(\alpha^{-n-1+k} u)
]+\beta^{-n}v
)
$ the corresponding truncation of $F^{-n}$; note that, if $\vert \alpha\vert>1$ and 
$\vert \beta\vert >1$,
then $F_N^{-n}$ is a polynomial contracting mapping.
\begin{theorem}
If $\vert \alpha\vert>1$, $\vert \beta\vert >1$ 
(thus $\{F^{\circ n}\}$ is {\sf not} a normal family) and $\vert\beta\vert<\vert\alpha\vert^N$,
then $\{F^{\circ n}\circ F_N^{-n}\}$ 
converges normally to a triangular
automorphism of the form
$G(u,v)=(u, \psi(u)+v)$ for a suitable entire function $\psi$ (hence $G(\CI^2)=\CI^2$).
\labelle{principal}
\end{theorem}
{\bf Proof:} trivially $\left(F^{\circ n}\circ F_N^{-n}   \right)_1(u,v)\equiv u$ and
\begin{eqnarray*}
\left(F^{\circ n}\circ F_N^{-n}\right)_2(u,v)
&=&
\sum_{k=0}^{n-1} \beta^k h(\alpha^{-1-k} u)-
\sum_{k=1}^n \beta^{n-k}P_N(\alpha{-n-1+k} u) +v\\
&=&
\sum_{k=0}^{n-1} \beta^k h(\alpha^{-1-k} u)-
\sum_{k=0}^{n-1}  \beta^{k}P_N(\alpha{-1-k} u) +v\\
&=&
\sum_{k=0}^{n-1}  \beta^{k}R_N(\alpha{-1-k} u) +v:=\psi_n(u)+v,
\end{eqnarray*}
where $R_N$ is the $N$-th remainder in the development of $h$.

Now 
$
\psi_n(u)=\sum_{k=0}^{n-1}\beta^k
\sum_{l=N}^{\infty}
\eta_l
(\alpha^{k+1})^{-l} u^l
%\labelle{stella}
$;
since, for $N\geq l$, we have
$$
\vert
{\beta^k}{\alpha^{-(k+1)l}}\vert
\leq
   \vert
{\beta}{\alpha^{-N}}
\vert^k
,$$
we can let $n$ diverge and exghange 
the order of summation, getting, uniformly on 
compact sets:
\begin{eqnarray}
\psi(u)&:=&
\lim_{n\to\infty}\psi_n(u) =
\sum_{l=N}^{\infty}
(\sum_{k=0}^{\infty}
 \beta^k (\alpha^{k+1})^{-l}
)
\eta_l u^l
\nonumber\\
&=&[
\sum_{k=0}^{\infty}
(\beta \alpha^{-l} )^k]\alpha^{-l}
\eta_l u^l=
\sum_{l=N}^{\infty}
(\alpha^l-\beta)^{-1}\eta_l u^l
;\labelle{stella}
\end{eqnarray}
by comparison with 
$h=\sum_{l=0}^{\infty}\eta_l u^l$, 
we see that the last series in (\ref{stella})
represents an entire function
$\psi$.
\QUAN
\vskip0,2truecm
{\bf An application:} let now $H$ be an automorphism of $\CI^2$, with a repulsive fixed point $p$,
$B$ the attracting basin of $p$ with respect to $H^{-1}$: by the theorem in the appendix of \cite{rudin}, there exists
a biholomorphic mapping $\Psi:B\rightarrow \CI^2$ and an elementary one $T:\CI^2\rightarrow\CI^2$ such that $H=\Psi^{-1}\CIRC T\CIRC\Psi$ on $B$: thus, if $\{f_n\}$ is a renormalizing polynomial family for $\{T^{\CIRC n}\}$ such that $\{T^{\CIRC n}\CIRC f_n\}$
admits the normal limit
$\Ef$, then, as in section 
\ref{mapm}, it is easily seen that $\{\psi^{-1}\CIRC f_n\}$ is a renormalizing family for $\{H^{\CIRC n}\}$
and $\{H^{\CIRC n}\CIRC \psi^{-1}\CIRC f_n\}$ 
converges normally to $\Psi^{-1}\CIRC\Ef$.
\vskip0.2truecm

{\bf Remark:} we end this section by recalling that,
if $f$ is a polynomial automorphism of $\CI^2$ whose iterates have uniformly bounded 
(algebraic)
degrees, then $f$ is conjugate to an elementary
mapping $F$, i.e. there exists a polynomial
automorphism $\phi$ of $\CI^2$ such that
$f=\phi^{-1}F\phi$ (see \cite{sibony}, pag. 123): then,
as a consequence of our preceeding theorem,
$f$ too admits a polynomial renormalizing family.
\section{\bf\large Elementary correspondences}
Following the same lines of reasoning,
theorem \ref{principal} could be extended to the case in which $h$ is the 'sum' of an entire function and a 
'Puisseux polynomial' with no 'branch points'
at $0$.
We shall also assume, as a simolifying hypothesis, that $h$
has a fixed point, as a corrspondence (see definition \ref{fixbra}) at $0$.
More precisely, let 
$\widehat{\h}$ be the germ (at $0$) of an entire function and
$\alpha=
\sum_{i=1}^M\sum_{l=1}^{N(i)}
\alpha_{il}(z-\zeta_i)^{\lambda_{li}}$,
as a correspondence, where $\{\zeta_i\}\subset\CI\setminus\{0\}$,
$\{\lambda_l\}\in\QU^+$
(note $\alpha$ is, according to
\cite{ahlfors}, pag. 300-306,
an 'algebraic function').
Now $\alpha$ posesses a branch
in a neighbourhood of $0$, identifying a 
holomorphic germ
$\alph$.
Let $\h=\widehat{\h} +\alph$  and $\G=\L+(0,\h)$, where $L$ is an expanding linear automorphism of $\CI^2$, with germ
$\L$; 
suppose that $\h(0)=0$ and let ${\cal R}$ be the elementary correspondence generated by $\G$
(see definition \ref{elemmap}). 
By construction, there exists exactly one branch,
 $\phi_n$ of ${\cal R}^{\CIRC n}$, defined on a neighbourhood $V_n$
of $0$,
 such that $\phi_n(0)=0$.
Then we have:
 \begin{theorem}
The family $\{\phi_n\}$ admits a renormalizing
family $\{\chi_n\}$, 
where $\chi_n=L^{-n}-(0,\alpha_n)$
for a suitable 
algebraic function $\alpha_n$;
the family $\{\phi_n\CIRC\chi_n\}$
converges normally to a triangular automorphism of the form
$G(u,v)=(u, \psi(u)+v)$ for a suitable entire function $\psi$ (hence $G(\CI^2)=\CI^2$).
\labelle{last}
\QUAN
\end{theorem}
{\bf Remark:} 
the proof proceeds almost verbatim as in theorem
\ref{principal}: 
the only difference is the degree of the polynomial truncations
of $h$:
we should indeed consider an integer 
$N$ so large that
$R_N$ contains no fractional-degree
terms.

\end{document}